\documentclass[12pt]{article}
\usepackage{amsmath,amssymb,amsthm,graphicx}
\numberwithin{equation}{section}

\newcommand\g{\gamma}
\renewcommand\d{\delta}

\newcommand\G{\Gamma}
\newcommand\f{\frac}
\newcommand{\Z}{{\mathbb{Z}}}
\newcommand{\R}{{\mathbb{R}}}

\newcommand{\A}{{\mathbb{A}}}

\newcommand{\U}{{\mathbb{H}}}

\newcommand{\quo}[1]{SL_{{#1}}({\Z})\backslash SL_{{#1}}({\R})/
SO_{{#1}}({\R})}

\renewcommand\i{^{-1}}

\renewcommand\({\left(}         % Quick delimiters
\renewcommand\){\right)}

\begin{document}

\begin{center}
{\Large A Summary of the Langlands-Shahidi Method of Constructing
L-functions}

\vspace{.7 cm}

{\sc Stephen D. Miller, Rutgers University}\footnote{Supported by
the National Science Foundation under Grant No. 0122799}

 {\tt
miller@math.rutgers.edu}

\vspace{.3 cm}
 December 2001
\end{center}

\vspace{.5 cm}

\begin{abstract}
These notes are from the Database of Automorphic L-functions at
{\tt http://www.math.rutgers.edu/$\sim$sdmiller/l-functions} .
They were written up by Stephen Miller, and are based on
discussions with Freydoon Shahidi, Purdue University.  They are
meant to serve as an introduction to the Langlands-Shahidi method
of studying L-functions through the Fourier coefficients of
Eisenstein series.

\end{abstract}

\section{Introduction}

The Langlands-Shahidi method uses information from spectral theory
on non-compact spaces to obtain the functional equations and many
analytic properties of L-functions, in particular several
important examples of Langlands L-functions.  The reason for the
connection to spectral theory is that L-functions arise in the
Fourier expansions of Eisenstein series.  The simplest example of
this phenomenon occurs already on $SL_2({\Z})\backslash \U$ for
the non-holomorphic Eisenstein series

$$E(x+iy,s)=\sum_{\stackrel{m,n\in\Z}{(m,n)=1}}
\f{y^s}{|m(x+iy)+n|^{2s}},$$

which has the Fourier expansion

\begin{equation}\label{eisfourexp}
  E(x+iy,s)=\sum_{n\in\Z} a_n(y,s)e^{2 \pi i n x},
\end{equation}
$$a_0(y,s)=y^s+\f{\xi(2s-1)}{\xi(2s)}y^{1-s},$$
$$a_n(y,s)=\f{2|n|^{s-1/2}\sigma_{1-2 s}(|n|)\sqrt{y}K_{s-1/2}(2
\pi |n|y)}{\xi(2s)},~n\neq 0,$$ where
$\xi(s)=\pi^{-s/2}\G(s/2)\zeta(s)=\xi(1-s)$ is the completed
Riemann $\zeta$ function, $\sigma_s(n)=\sum_{d|n}d^s,$ and
$K_s(y)$ is the K-Bessel function. One can deduce from the
functional equation $\xi(s)=\xi(1-s)$ that

\begin{equation}\label{eisfunc}
  E(z,s)=\f{\xi(2s-1)}{\xi(2s)}E(z,1-s).
\end{equation}

The main point of using spectral theory here is {\em Maass'
Lemma}, which ultimately implies that because $E(z,s)$ and
$E(z,1-s)$ both have Laplace eigenvalue $s(1-s)$,  the two must be
multiples of each other.  The ratio can be found from the constant
term $a_0(y,s)$, and so the functional equation (\ref{eisfunc})
can be proven without knowing the functional equation
$\xi(s)=\xi(1-s)$.

 Many results from functional analysis have been
used to generalize the analytic continuation and functional
equation of Eisenstein series, starting with Selberg and continued
by Langlands \cite{langeis}.  Their spectral methods in general do
not assume any information about the L-functions (generalizing
$\zeta(s)$) that are involved.  In his Yale Monograph {\it Euler
Products} (\cite{langeul}), Langlands obtained the meromorphic
continuation of a wide variety of L-functions using the theory of
the constant term.  His calculations led to the definition of the
L-group, and ultimately to the formulation of his functorality
conjectures.

Furthermore, the functional equation (\ref{eisfunc}) can be turned
around to gain information about the symmetries of the terms in
the Fourier expansion (\ref{eisfourexp}).  In particular, an
analysis of the first Fourier coefficient readily gives the
functional equation $\xi(s)=\xi(1-s)$.  This latter observation
has been developed by Shahidi; in conjunction with Langlands'
constant term theory, he has applied an analysis to the
non-constant terms of general Eisenstein series to obtain the
functional equations and meromorphic continuation (with only a
finite number of poles between $1/2$ and 1) of the L-functions
which occur in the Fourier expansions of Eisenstein series.

In general it has been a difficult challenge to prove the
L-functions are entire.  Kim has applied the following idea from
representation theory to rule out the poles on the real axis
between $1/2$ and 1 in many cases.  These potential poles are also
singularities of Eisenstein series, and their residues are
$L^2$-automorphic forms, which always correspond to unitary
representations that can be explicitly described by the Eisenstein
series they came from.  Kim remarked that known results about the
unitary dual show that many of these representations do not occur,
allowing one to conclude the holomorphy of the Eisenstein series
and the L-functions at these points!  When combined with
\cite{gelsha}, \cite{shaduke}, and \cite{shannals}, this has
recently led to new examples of entire L-functions with
breakthrough applications to the Langlands functorality
conjectures (\cite{CKPSS},\cite{K-S}).

\section{An Outline of the Method}

The following is a brief sketch of the main points of the method;
a fuller introduction with more definitions and detailed examples
can be found in \cite{korea}.  Detailed examples of constant term
calculations can be found in many places, e.g.
\cite{langeul},\cite{bluelands},\cite{gelsha}, and \cite{mil}.
Though it is possible to describe the method without adeles (as
was done in the introduction), their use is key in higher rank for
factoring infinite product expansions into L-functions.

Let $F$ be a field, $\A=\A_F$ its ring of adeles, and $G$ a split
group over $F$.  Much carries over to quasi-split case as well,
and we will highlight the technical changes needed for this at the
end.  Fix a Borel (= a maximal connected solvable) subgroup
$B\subset G$, and a standard maximal parabolic $P\supset B$
defined over $F$.\footnote{The theory can be extended to
non-maximal parabolics, but this does not yield any extra
information about L-functions.}  Decompose $B=TU,$ where $T$ is a
maximal torus.  The parabolic can be also decomposed as $P=MN$,
where the unipotent radical $N\subset U$, and
 $M$ is the unique Levi component containing $T$.  Denote by
$^LG,^LM,^LN,$ etc. the Langlands dual L-groups (see
\cite{korea}).

One of the key aspects of this method is that it uses many
possibilities of parabolics of different groups, especially
exceptional groups.  This is simultaneously a strength (in that
there is a wide range of possibilities) and a limitation (in that
there are only finitely many exceptional groups).

\subsection{Cuspidal Eisenstein Series}

Recall that an automorphic form in $L^2(\Gamma\backslash G)$ is
associated to a (unitary) automorphic representation of $G$
mapping $\phi(h)\mapsto \phi(hg)$.  Let $\pi=\otimes_v\pi_v$ be a
cuspidal automorphic representation of $M(\A)$; we may assume that
almost all components $\pi_v$ are spherical unitary
representations (meaning that they have a vector fixed by
$G(O_v)$, where $O_v$ is the ring of integers of the local field
$F_v$).  For these places $v$ the equivalence class of the unitary
representation $\pi_v$ is determined by a semisimple conjugacy
class $t_v\in~^LG$, the L-group.  This conjugacy class is used to
define the L-functions below in (\ref{lfuncdef}).

A maximal parabolic has a modulus character $\d_P$, which is the
ratio of the Haar measures on $M\cdot N$ and $N\cdot M$.  It is
related to the simple root of $G$ which does not identically
vanish on $P$.  For any automorphic form $\phi$ in the
representation space of $\pi$, we can define the Eisenstein series

\begin{equation}\label{eisdef}
  E(\pi,s,g)=\sum_{\g\in P(F)\backslash G(F)} \phi(\g g)\d_P(\g
  g)^s,
\end{equation}
and their constant terms
\begin{equation}\label{constdef}
 c(\pi,s,g) =  \int_{N'(F)\backslash N'(A)} E(\phi,s,ng)dn,
\end{equation}
where $N'$ is the {\em opposite} parabolic to N (it is related by
the longest element in the Weyl group).

\subsection{Langlands L-functions}

If $\rho$ is a finite-dimensional complex representation of $^LM$
and $S$ is a finite set including the archimedean and ramified
places, then the partial Langlands L-function is

\begin{equation}\label{lfuncdef}
  L_S(s,\pi,\rho)=\prod_{v \notin S} \det(I-\rho(t_v)q_v^{-s})\i.
\end{equation}
Here $q_v$ is the cardinality of the residue field of $F_v$, a
prime power.  The full, {\em completed}, L-function involves extra
factors for places in $S$, whose definition is technical and in
general difficult.

\subsection{The Constant Term Formula}

The constant term formula involves the sum of two terms.  The
first, which only occurs when the parabolic is its own opposite,
is essentially just $\phi(g)$.  Langlands showed that the map to
the second term is given by an operator
\begin{equation}\label{consttermformula}
M(s,\pi)=\(\prod_{j=1}^m
\f{L(a_js,\tilde{\pi},r_j)}{L(1+a_js,\tilde{\pi},r_j)}\)\otimes_{v\in
S} A(s,\pi_v),\end{equation} where $A(s,\pi_v)$ are a finite
collection of operators, $r$ the adjoint action of $^LM$ on the
lie algebra of $^LN$, $r_1,\ldots,r_m$ the irreducible
representations it decomposes into, and $a_j$ integers which are
multiples of each other (coming from roots related to the $r_j$).
See \cite{korea} for a fuller discussion along with an example for
the Lie group $G_2$ and the symmetric cube L-function. Tables
listing Lie groups and the representations $r_j$ occurring for
them can be found in \cite{langeul} and \cite{shannals}.

\subsection{The Non-Constant Term: Local Coefficients}

We must now make a further restriction on the choice of $\pi$
involved, namely that it be {\em generic}, i.e. have a Whittaker
model.  This means that if $\psi$ is a generic unitary character
of $U(F)\backslash U(A)$, we need to require

$$W_v(g,\psi)=\int_{U_M(F)\backslash U_M(\A)}\phi(ng)\overline{\psi(n)}dn \neq 0   ,
~~~U_M=U\cap M$$
for some $\phi$ and $g$.

Shahidi's formula uses the Casselman-Shalika formula for Whittaker
functions to express the non-constant term as

\begin{equation}\label{nonconstform}
  \int_{N'(F)\backslash N'(A)}E(ng)\overline{\psi(n)}dn = \prod_{v\in
  S}W_v(1)\prod_{j=1}^m\f{1}{L(1+a_js,\tilde{\pi},r_j)}.
\end{equation}

Applying the functional equation of the Eisenstein series (which
has the constant-term ratio involved), one gets the ``crude''
functional equation for the product of $m$ L-functions

\begin{equation}\label{bigprod}
  \prod_{j=1}^m L_S(a_j
  s,\tilde{\pi},r_j)=\prod_{j=1}^mL_S(1-a_js,\pi,r_j)\prod_{v\in S}\(\mbox{local
  factors}\).
\end{equation}

Shahidi's papers \cite{shaduke} and \cite{shannals} match all the
local factors above to the desired L-functions (as in the remark
after (\ref{lfuncdef})). This gives the full functional equation
for these $m$ L-functions, but only when multiplied together. His
1990 paper \cite{shannals} uses an induction to isolate each of
the $m$ factors above separately.

\subsection{Analytic Properties and the Quasi-Split Case}

It remains to prove that the L-functions are entire, except
perhaps at $s=0$ and 1 where the order of the poles is understood.
The theory of Eisenstein series provides this full analyticity
except for when $\pi$ satisfies a self-duality condition; even in
this case, it can be shown that the L-functions have only a finite
number of poles, all lying on the real axis between $\f 12$ and 1.
 Kim's observation of using the unitary dual has
worked in many cases. It is also always possible to remove the
potential poles by twisting by a highly-ramified $GL(1)$ character
of $\A_F$; this has been crucial for applications to functorality
by the converse theorem \cite{converse}.

The main difference in the quasi-split case is that the action of
the Galois group $G_F$ is no longer trivial.   The L-groups are
potentially disconnected as a semi-direct product of a connected
component and $G_F$.  Also, the representation $\rho$ used to
define the Langlands L-functions may also depend on the place $v$.


\begin{thebibliography}{99}

\bibitem[CKPSS]{CKPSS} {\sc J. W. Cogdell, H.H. Kim,
I.I. Piatetski-Shapiro} and {\sc
 F.   Shahidi}, {\it On lifting from classical groups to ${\rm GL}\sb
 N$},
   Inst. Hautes Études Sci. Publ. Math. {\bf 93} (2001), 5--30.

\bibitem[CPS]{converse} {\sc J.W.
Cogdell} and {\sc I.I. Piatetski-Shapiro}, {\it Converse
   theorems for ${\rm GL}\sb n$, II}, J. Reine Angew. Math. {\bf 507}
    (1999),
   165--188.

\bibitem[G-S]{gelsha} {\sc Stephen Gelbart} and {\sc Freydoon Shahidi},
{\it Analytic
   properties of automorphic $L$-functions},  Perspectives in Mathematics,
   {\bf 6}, Academic Press, Inc., Boston, MA, 1988. viii+131 pp.

\bibitem[K-S]{K-S} {\sc Henry H. Kim} and {\sc Freydoon Shahidi},
{\it Functorial products
   for $\rm GL\sb 2\times GL\sb 3$ and functorial symmetric cube for $\rm
   GL\sb 2$}, C. R. Acad. Sci. Paris Sér. I Math. {\bf 331} (2000), no. 8,
   599--604.

\bibitem[L1]{langeis} {\sc Robert P. Langlands}, {\it On the
Functional Equations Satisfied by Eisenstein Series}, Springer
Lecture Notes in Mathematics {\bf 544}, 1976.
 \newline
 Available at {\tt
http://sunsite.ubc.ca/DigitalMathArchive/Langlands/}

\bibitem[L2]{langeul} Robert P. Langlands, {\em Euler Products},
Yale Mathematical Mongraphs, {\bf 1}, 1971.
 \newline
 Available at {\tt
http://sunsite.ubc.ca/DigitalMathArchive/Langlands/}

\bibitem[L3]{bluelands} {\sc R. P. Langlands}, {\it Eisenstein series, the trace
   formula, and the modern theory of automorphic forms}, Number theory,
   trace formulas and discrete groups (Oslo, 1987), 125--155, Academic
   Press, Boston, MA, 1989.

\bibitem[M]{mil} {\sc Stephen D. Miller}, {\it Cusp Forms on
$\quo{3}$},
 and lecture notes with
calculations, to be posted on the database website.  {\tt
http://www.math.rutgers.edu/$\sim$sdmiller/l-functions}

\bibitem[Sh1]{korea} {\sc Freydoon Shahidi}, {\it Intertwining Operators, L-functions, and
Representation Theory}, Lecture Notes of the Elevent KAIST
Mathematics Workshop (1996), Ja Kyung Koo, ed., pp. 1-63.

\bibitem[Sh2]{shaduke} {\sc Freydoon Shahidi}, {\it Local coefficients as Artin
 factors for real groups}, Duke Math. J. {\bf 52} (1985), no. 4, 973--1007.

\bibitem[Sh3]{shannals} {\sc Freydoon Shahidi}, {\it
   A proof of Langlands' conjecture on Plancherel measures; complementary
   series for $p$-adic groups},
   Ann. of Math. (2) {\bf 132} (1990), no. 2, 273--330.



\end{thebibliography}
\end{document}